\documentclass[11pt]{article}
\renewcommand{\baselinestretch}{1.2}
\newcommand{\dated}{\mbox{} \hfill {\small [{\tt \today}]}} % theorems etc.
\usepackage{amsthm,enumerate}
\theoremstyle{plain}
\newtheorem{theorem}{Theorem}[section]
\newtheorem{lemma}[theorem]{Lemma}
\newtheorem{corollary}[theorem]{Corollary}
\newtheorem{proposition}[theorem]{Proposition}
\theoremstyle{definition}
\newtheorem{definition}[theorem]{Definition}
\theoremstyle{remark}
\newtheorem*{remark}{Remark}
\newtheorem*{example}{Example}
\newtheorem*{rems}{Remarks}
\newtheorem*{exs}{Examples}
\newenvironment{remarks}{\begin{rems}\begin{enumerate}}{\end{enumerate}\end{rems}}

\newenvironment{items}{\begin{enumerate}[\rm (i)]}{\end{enumerate}}
\newenvironment{alphitems}{\begin{enumerate}[\rm (a)]}{\end{enumerate}}

 \usepackage{amsmath,amssymb,amsfonts}
\newenvironment{keywords}{\noindent\small {\it Keywords\/}:}{\vskip 4pt}
\newenvironment{classification}{\noindent\small 2000 {\it Mathematics Subject
Classification\/}:}{\vskip 12pt}

\newcommand{\defiff}{\quad:\Longleftrightarrow\quad}
\newcommand{\comps}{{\mathbb C}}
\newcommand{\reals}{{\mathbb R}}

\newcommand{\free}{{\mathbb F}}

\newcommand{\void}{\varnothing}
\newcommand{\tensor}{\otimes}

\newcommand{\cstar}{{C^\ast}}

\newcommand{\cb}{{\mathrm{cb}}}

\newcommand{\A}{{\mathfrak A}}

\newcommand{\Hilbert}{{\mathfrak H}}

\newcommand{\varcl}[1]{\overline{#1}}

\newcommand{\AM}{\mathrm{AM}}
\newcommand{\Cb}{\mathcal{C}_b}
\newcommand{\Mcb}{M_{cb}}
\newcommand{\Mcbo}{M_{cb,0}}
\newcommand{\AMcb}{A_{M_{cb}}}
\newcommand{\AP}{\mathcal{AP}}
\newcommand{\WAP}{\mathcal{WAP}}
\newcommand{\SO}{\mathrm{SO}}
\title{(Non-)amenability of the Fourier algebra \\ in the $cb$-multiplier norm}
\author{\textit{Volker Runde}\thanks{Research supported by NSERC.}}
\date{}
\begin{document}
\maketitle
\begin{abstract}
For a locally compact group $G$, let $A(G)$ denote its Fourier algebra, $\Mcb(A(G))$ the completely bounded multipliers of $A(G)$, and $\AMcb(G)$ the closure of $A(G)$ in $\Mcb(A(G))$. We show that, if $\AMcb(G)$ is amenable, then $a(G_d)$, the almost periodic compactification of the discretization of $G$, has an abelian subgroup of finite index. As a consequence, $\AMcb(G)$ cannot be amenable if $G$ contains a copy of $\free_2$, the free group in two generators, as a closed subgroup.
\end{abstract}
\section*{Introduction}
The Fourier algebra $A(G)$ of a locally compact group $G$ was introduced by P.\ Eymard in \cite{Eym}; if $G$ is abelian with dual group $\hat{G}$, these algebras are isometrically isomorphic to $L^1(\hat{G})$, the group algebra of $\hat{G}$. In \cite{Lep}, H.\ Leptin characterized the amenable, locally compact groups in terms of their Fourier algebras: a locally compact group $G$ is amenable if and only if $A(G)$ has a bounded approximate identity (in fact, of bound one).
\par 
In \cite{Joh1}, B.\ E.\ Johnson initiated the theory of amenable Banach algebras and proved that a locally compact group $G$ is amenable if and only if its group algebra $L^1(G)$ is amenable (\cite[Theorem 2.5]{Joh1}). As all amenable Banach algebras have bounded approximate identities, this prompted the conjecture that $A(G)$ is amenable if and only if $G$ is amenable. Alas, in \cite{Joh3}, Johnson showed that there are compact groups $G$---among them $\SO(3)$---for which $A(G)$ fails to be amenable. In \cite{FR1} (see also \cite{RunPAMS}), B.\ E.\ Forrest and the author finally showed that $A(G)$ is amenable if and only if $G$ is almost abelian, i.e., has an abelian subgroup of finite index.
\par 
It turns out that, if one wants to obtain a satisfactory amenability theory for Fourier algebras, one has to take their canonical operator space structure into account. Such a theory was initiated in \cite{Rua} by Z.-J.\ Ruan under the name \emph{operator amenability}. He was able to show: a locally compact group $G$ is amenable if and only if $A(G)$ is operator amenable (\cite[Theorem 2.5]{Rua}).
\par 
The operator space structure of the Fourier algebra $A(G)$ of a locally compact group $G$, can be used to define another norm on it. A multiplier on $A(G)$ is a function $\phi \!: G \to \comps$ such that $\phi \!: G \to \comps$ such that $\phi A(G) \subset A(G)$. We consider the $cb$-multipliers on $A(G)$, i.e., those multipliers $\phi$ such that the map $A(G) \ni f \mapsto \phi f$ is completely bounded. They form a closed, commutative subalgebra of the completely bounded operator on $A(G)$, which we denote by $\Mcb(A(G))$. As $A(G)$ embeds (completely) contractively into $\Mcb(A(G))$, it inherits a new norm from $\Mcb(A(G))$, the \emph{$cb$-multiplier norm}. For amenable $G$, this norm equals the given norm (as a consequence of \cite{Lep}). For non-amenable groups, however, the two norms are inequivalent. We denote the closure of $A(G)$ in $\Mcb(A(G))$ by $\AMcb(G)$. Unlike $A(G)$, the algebra $\AMcb(G)$ may have a bounded approximate identity for non-amenable $G$: this is the case, for instance, if $G = \free_2$, the free group in two generators (\cite{dCH}). Even more surprisingly, $\AMcb(\free_2)$ is operator amenable as shown by B.\ E.\ Forrest, N.\ Spronk, and the author (\cite{FRS}). A crucial ingredient of the proof was the fact $\Mcb(A(\free_2))$ contains a weak$^\ast$ dense subalgebra that is completely isometrically isomorphic to the operator amenable completely contractive Banach algebra $B(a\free_2)$ with $a\free_2$ denoting the almost periodic compactification of $\free_2$. 
\par 
The result from \cite{FRS} naturally begets the question if $\AMcb(\free_2)$ is amenable in the sense of Johnson's original definition. The methods used in \cite{FRS} are not suited to produce an affirmative answer to this question. On the negative side, Forrest and the author were able to show that a locally compact group $G$ has to be abelian if $\AM(\AMcb(G)) < \frac{2}{\sqrt{3}}$, where $\AM(\AMcb(G))$ denotes the amenability constant of $\AMcb(G)$ as introduced in \cite{Joh3}. 
\par 
In the present paper, we prove that the amenability of $\AMcb(G)$ forces $a(G_d)$, i.e., the almost periodic compactification of the discretization of $G$, to be almost abelian. As in \cite{FRS}, almost periodic $cb$-multipliers play a pivotal r\^ole in the proof. If $G$ is discrete and maximally almost periodic, but fails to be almost abelian, then $\AMcb(G)$ cannot be amenable; in particular, $\AMcb(\free_2)$ is not amenable. As a consequence, $\AMcb(G)$ is not amenable if it contains a closed subgroup isomorphic to $\free_2$.
\par 
Even though the main result of the paper is entirely Banach algebraic, its object is defined in terms of operator space theory. It should therefore come as no surprise that---as in \cite{FR1} and \cite{RunPAMS}---operator spaces play a crucial r\^ole in its proof. 
\section{Preliminaries}
Amenable Banach algebras were introduced by B.\ E.\ Johnson in \cite{Joh2}. The motivation for the choice of terminology is \cite[Theorem 2.5]{Joh1}: a locally compact group $G$ is amenable if and only if its group algebra $L^1(G)$ is amenable. We will not work with the definition from \cite{Joh1}, but with a characterization of amenable Banach algebras from \cite{Joh2} instead.
\par
Following \cite{ER}, we use $\tensor^\gamma$ to denote the completed projected tensor product of Banach spaces. Given a Banach algebra $\A$, a left Banach $\A$-module $E$, and a right Banach $\A$-module $F$, the tensor product $E \tensor^\gamma F$ becomes a Banach $\A$-bimodule via
\[
  a \cdot (x \tensor y) := a \cdot x \tensor y \quad\text{and}\quad (x \tensor y) \cdot a := x \tensor y \cdot a \qquad (a \in \A, \, x \in E, \, y \in F).
\]
In particular, $\A \tensor^\gamma \A$ is a Banach $\A$-bimodule in a canonical manner.
\par 
The following is from \cite{Joh2}:
\begin{theorem} \label{amthm}
A Banach algebra $\A$ is amenable if and only if it has a \emph{bounded approximate diagonal}, i.e., a bounded net $( \boldsymbol{d}_\alpha )_\alpha$ in $\A \tensor^\gamma \A$ such that
\[
  a \cdot \boldsymbol{d}_\alpha - \boldsymbol{d}_\alpha \cdot a \to 0 \qquad (a \in \A)
\]
and 
\[
  a \Delta \boldsymbol{d}_\alpha \to a \qquad (a \in \A)
\]
where $\Delta \!: \A \tensor^\gamma \A \to A$ is the multiplication map given by
\[
  \Delta(a \tensor b) = ab \qquad (a,b \in \A).
\]
\end{theorem}
\par 
Our main reference for the theory of amenable Banach algebras is \cite{RunAM}.
\begin{remarks}
\item If $\A$ has an approximate diagonal bounded by $C \geq 1$, we call $\A$ \emph{$C$-amenable}.
\item The \emph{amenability constant} $\AM(\A)$ of $\A$ is defined as
\begin{equation} \label{amconst} 
  \AM(A) := \inf \{ C \geq 1 : \text{$\A$ is $C$-amenable} \}.
\end{equation}
\item It is easy to see that the infimum (\ref{amconst}) is attained, so that every amenable Banach algebra $\A$ is $\AM(A)$-amenable.
\end{remarks}
\par
Let $G$ be a locally compact group. The Fourier algebra $A(G)$ of $G$ was defined by P.\ Eymard in \cite{Eym}. A more recent, comprehensive monograph is \cite{KL}.
\par
A \emph{multiplier} of $A(G)$ is a function $\phi \!: G \to \comps$ such that $\phi A(G) \subset A(G)$. A straightforward application of the closed graph theorem yields that, if $\phi \!: G \to \comps$ is a multiplier of $A(G)$, then the operator
\[
  M_\phi \!: A(G) \to A(G), \quad f \mapsto \phi f
\]
is bounded. As $A(G)$ is an operator space in a natural manner (\cite[Section 16.2]{ER}), it makes sense to define:
\begin{definition}
Let $G$ be a locally compact group. Then a multiplier $\phi \!: G \to \comps$ is called \emph{completely bounded} if $M_\phi \!: A(G) \to A(G)$ is completely bounded. We denote the collection of completely bounded multipliers---short: \emph{$cb$-multipliers}---of $A(G)$ by $\Mcb(A(G))$. We set
\[
  \| \phi \|_{\Mcb} := \| M_\phi \|_{cb} \qquad (\phi \in \Mcb(A(G))).
\]
\end{definition} 
\par
Here, $\| \cdot \|_{cb}$ stands for the completely bounded norm (\cite{ER}).
\begin{remarks}
\item It is clear that $(\Mcb(A(G)), \| \cdot \|_{\Mcb})$ is a completely contractive---see \cite[Section 16.1]{ER}---Banach algebra.
\item In \cite{dCH}, it was observed that $\Mcb(A(G))$ is a dual Banach space in a canonical way for every locally compact group $G$. The resulting weak$^\ast$ topology coincides on norm bounded subsets of $\Mcb(A(G))$ with the relative topology induced by $\sigma(L^\infty(G), L^1(G))$; by the Kre\u{\i}n--\v{S}mulian Theorem (\cite[Theorem V.6.4]{DS}), this property uniquely identifies this particular predual. It is straightforward to see that this duality turns $\Mcb(A(G))$ into a dual Banach algebra in the sense of \cite{RunSt}. (In \cite{Spr}, N.\ Spronk provided an alternative characterization of this canonical predual.) 
\end{remarks}
Let $G$ be a locally compact group. As $A(G)$ is a completely contractive Banach algebra, it is clear that $A(G)$ embeds (completely) contractively into $\Mcb(A(G))$. We thus define:
\begin{definition}
Let $G$ be a locally compact group. Then $\AMcb(G)$ is defined to be the norm closure of $A(G)$ in $\Mcb(A(G))$. 
\end{definition}
\begin{remarks}
\item It follows from Leptin's Theorem (\cite{Lep}) that $\| \cdot \|_{\Mcb}$ and the given norm on $A(G)$ are identical for amenable $G$. Conversely, if $\| \cdot \|_{cb}$ and the given norm are equivalent on $A(G)$, then $G$ is amenable: this seems to have been known for quite some time (see  the remarks in \cite{Spr} following Proposition 6.9), However, a correct proof accessible to the public came out only as late as in \cite{LX}.
\item In \cite{dCH} it is shown that $A(\free_2)$ has an approximate identity that is bounded with respect to $\| \cdot \|_{\Mcb}$, so that $\AMcb(\free_2)$ has a bounded approximate identity.
\end{remarks}
\par 
Given a group $G$, we call the set 
\[
  \nabla(G) := \{ (x,x^{-1}) : x \in G \}
\]
the \emph{anti-diagonal} of $G$. For any subset $S$ of $G$, we denote its indicator function by $\chi_S$. 
\par
The following result, which will be crucial for our argument, is \cite[Proposition 3.3]{FR2}:
\begin{proposition} \label{BrianVolker}
Let $G$ be a locally compact group such that $\AMcb(G)$ is amenable. Then we have $\chi_{\nabla(G)} \in \Mcb(A(G_d \times G_d))$ with $\| \chi_{\nabla(G)} \| \leq \AM(\AMcb(G))$.
\end{proposition}
\par 
Here, $G_d$ stands for the discretization of $G$, i.e., $G$ equipped with the discrete topology.
\section{$wG$ and $aG$ for a locally compact group $G$}
The following definitions are standard:
\begin{definition}
A \emph{semitopological semigroup} is a semigroup $S$ equipped with a Hausdorff topology if the multiplication
\[
  S \times S \to S, \quad (s,t) \mapsto st
\]
is separately continuous. 
\end{definition}
\begin{remarks}
\item If multiplication in $S$ is jointly continuous, we call $S$ a \emph{topological semigroup}.
\item If $S$ has an identity, we denote it by $e_S$.
\end{remarks}
\begin{definition}
Let $S$ be a semitopological semigroup. An \emph{involution} on $S$ is a continuous map $^\ast \!: S \to S$ such that
\[
  s^{\ast\ast} = s \qquad (s,t \in S)
\]
and
\[
  (st)^\ast = t^\ast s^\ast \qquad (s,t \in S).
\]
\end{definition}
\par 
Given a topological semigroup $S$, a function $f \!: S \to \comps$, and $s \in S$, we define the left and right translates of $f$ by $s$ via
\[
  L_sf(t) := f(st) \quad\text{and}\quad R_s f(t) := f(ts) \qquad (t \in S).
\]
We denote the bounded continuous functions on $S$ by $\Cb(S)$---if $S$ is compact, we will simply write $\mathcal{C}(S)$---, and define $f \in \Cb(S)$ to be \emph{weakly almost periodic} if the set 
$\{ L_s f : s \in S \}$ is relatively weakly compact in $\Cb(S)$. (By Grothendieck's Double Limit Criterion, we could have defined weakly almost periodic functions in terms of right translates as well.)
\begin{definition} \label{wapdef}
Let $S$ be a semitopological semigroup. We set 
\[
  \WAP(S) := \{ f \in \Cb(S) : \text{$f$ is weakly almost periodic} \}.
\]
\end{definition}
\begin{remarks}
\item $\WAP(S)$ is a commutative $\cstar$-algebra with identity; we denote its compact character space by $wS$, the \emph{weakly almost periodic compactification} of $S$. There is a canonical continuous map $\iota_w \!: S \to wS$ with dense range. Moreover, $wS$ can be turned into a semitopological semigroup such that $\iota_w$ becomes a homomorphism of semigroups. See \cite{Bur} for this and other background information.
\item If $S$ is locally compact, then $\iota_w \!: S \to wS$ is a homeomorphism onto its image (\cite[Theorem 3.6]{Bur}).
\item If $G$ is a locally compact group, $wG$ need not be a group anymore. However, inversion in the group extends to $wG$ as an involution.
\item If $G$ and $H$ are locally compact groups, then $w(G \times H) = wG \times wH$ need not hold (\cite[Theorem 5.1]{Chou}). However, the identity on $G \times H$ extends to a continuous semigroup homomorphism from $w(G \times H)$ onto $wG \times wH$.
\end{remarks}
\par 
The following was stated in \cite{Xu}, where it was remarked that it followed from Gro\-then\-dieck's Double Limit Criterion. For the reader's convenience, we present a complete proof.
\begin{proposition} \label{McbWAPprop1}
Let $G$ be a locally compact group. Then $\Mcb(A(G)) \subset \WAP(G)$ holds.
\end{proposition}
\begin{proof}
Let $\phi \in \Mcb(G)$, and let $(x_n)_{n=1}^\infty$ and $(y_n)_{n=1}^\infty$ be sequences in $G$ such that the double limits $\lim_{n \to \infty} \lim_{m \to \infty} \phi(x_n y_m)$ and $\lim_{m \to \infty} \lim_{n \to \infty} \phi(x_n y_m)$ both exist. We need to show that they coincide.
\par
By \cite{Gil}---for a short and accessible proof, see \cite{Jol}---, there are a Hilbert space $\Hilbert$ and bounded, continuous functions $\boldsymbol{\xi}, \boldsymbol{\eta} \!: G \to \Hilbert$ such that
\[
  \phi(xy^{-1}) = \langle \boldsymbol{\xi}(x) | \boldsymbol{\eta}(y) \rangle \qquad (x,y \in G).
\]
(We use $\langle \cdot | \cdot \rangle$ for Hilbert space inner products as opposed to $\langle \cdot, \cdot \rangle$, which we use for Banach space duality.) The sequences $( \boldsymbol{\xi}(x_n) )_{n=1}^\infty$ and $(\boldsymbol{\eta}(y_n^{-1}) )_{n=1}^\infty$ are bounded in $\Hilbert$. As closed balls in $\Hilbert$ are weakly compact, there are $\xi, \eta \in \Hilbert$ as well as subnets $( x_{n_\alpha} )_\alpha$ and $( y_{n_\beta} )_\beta$ of $(x_n)_{n=1}^\infty$ and $(y_n)_{n=1}^\infty$, respectively, such that 
\[
  \xi = \text{weak-}\lim_\alpha \boldsymbol{\xi}(x_{n_\alpha}) \qquad\text{and}\qquad \eta = \text{weak-}\lim_\beta \boldsymbol{\eta}(y^{-1}_{n_\beta}).
\]
It follows that 
\[
  \lim_{n \to \infty} \lim_{m \to \infty} \phi(x_n y_m) = \lim_{n \to \infty} \lim_{m \to \infty} \langle \boldsymbol{\xi}(x_n) | \boldsymbol{\eta}(y_m^{-1}) \rangle 
  = \lim_\alpha \lim_\beta \langle \boldsymbol{\xi}(x_{n_\alpha}) | \boldsymbol{\eta}(y_{m_\beta}^{-1}) \rangle 
  = \langle \xi, \eta \rangle
\]
and---similarly---
\[
  \lim_{m \to \infty} \lim_{n \to \infty} \phi(x_n y_m) = \langle \xi, \eta \rangle,
\]
which completes the proof.
\end{proof}
\begin{remark}
In \cite{Xu}, the claim is stated for Herz-Schur multipliers in an $L^p$-context for $p \in (1,\infty)$. The same argument works there, too: the Hilbert spaces have be replaced with so-called $\mathit{QSL}_p$-spaces, i.e., quotients of subspaces of $L^p$-spaces. 
\end{remark}
\par 
Given a topological semigroup $S$, we define $f \in \Cb(S)$ to be \emph{almost periodic} if $\{ L_s f : s \in S \}$ is relatively (norm) compact in $\Cb(S)$. We define:
\begin{definition} 
Let $S$ be a semitopological semigroup. We set 
\[
  \AP(S) := \{ f \in \Cb(S) : \text{$f$ is almost periodic} \}.
\]
\end{definition}
\begin{remarks}
\item Like in the case of weak almost periodicity, $\AP(S)$ can as well be defined through right translates.
\item Like $\WAP(S)$, $\AP(S)$ is a commutative $\cstar$-algebra with identity; we denote its compact character space by $aS$, the \emph{almost periodic compactification} of $S$. There is a canonical continuous map $\iota_a \!: S \to aS$ with dense range. Moreover, $aS$ can be turned into a topological semigroup such that $\iota_a$ becomes a homomorphism of semigroups. See \cite[Note 1.14]{Bur} for this and other background information.
\item Through restriction, we have a canonical continuous semigroup homomorphism from $aS$ onto $aS$, which is a quotient map.
\item If $G$ is a locally compact group, then $aG$ is a compact group turning $\iota_a$ into a group homomorphism (\cite[Corollary 2.27]{Bur}). Note, however that $\iota_a \!: G \to aG$ need not be injective, e.g., if $G = \mathrm{SL}(2,\reals)$ (\cite[(22.22)(h)]{HR1}). Groups for which $\iota_a$ is injective are known as \emph{maximally almost periodic} (see \cite[3.2.17]{Pal}). Compact and abelian groups are maximally almost periodic, but so is $\free_2$, the free group in two generators (this follows from \cite[(2.9) Theorem]{HR1}). For abelian $G$, $aG$ is the usual Bohr compactification of $G$.
\item Unlike the weakly almost periodic compactification, the almost periodic compactification respects Cartesian products: If $G$ and $H$ are locally compact groups, then there is a canonical homeomorphic group isomorphism $a(G \times H) \cong aG \times aH$ (\cite[2.4 Theorem]{BJM}). 
\item For any locally compact group $G$, the semigroup $wG$ has a minimal ideal $\ker wG$, the \emph{kernel} of $wG$, which is a compact group (\cite[Theorem 2.7]{Bur}). Through restriction onto $\AP(G)$ it is canonically  homeomorphically isomorphic to $aG$ (\cite[Theorem 2.26]{Bur}). Let $e_{\ker wG}$ be the identity of $\ker wG$, and define
\[
  p_{\ker wG} \!: wG \to \ker wG, \quad s \mapsto e_{\ker wG}\, s.
\]
We obtain the following sequence of semigroup homomorphism
\[
  G \stackrel{\iota_a}{\rightarrow} aG \cong \ker wG \hookrightarrow wG \stackrel{p_{\ker wG}}{\twoheadrightarrow} \ker wG \cong aG,
\]
of which the last two arrows are the identity on $\ker wG$.
\item Let $s \in wG$. As $\ker wG$ is an ideal of $wG$, $s \,e_S$ and $e_S \, s$ lie in $\ker wG$, so that 
\[
  s \, e_{\ker wG} = e_{\ker wG} \, s \, e_{\ker wG} = e_{\ker wG} \, s \qquad (s \in wG), 
\]
i.e, $e_{\ker wG}$ is central in $wG$.
\end{remarks}
\section{$\Mcb(A(G))$ and $B(aG)$}
Let $G$ be a locally compact group. For $f \in \WAP(G)$, we denote by $\hat{f}$ its Gelfand transform on $wG$. Let
\[
  \WAP_0(G) := \{ f \in \WAP(G) : \hat{f} |_{\ker wG} \equiv 0 \}.
\]
Then $\WAP_0(G)$ is a closed ideal of $\WAP(G)$, and we have a direct sum decomposition
\[
  \WAP(G) = \AP(G) \oplus \WAP_0(G)
\]
(\cite[Theorem 2.22]{Bur}).
\par
Through composition, $\iota_a \!: G \to aG$ induces a (completely) isometric algebra homomorphism $\iota_a^\ast \!: B(aG) \to \Mcb(A(G))$ (\cite[Corollary 6.3(ii)]{Spr}). We can therefore canonically identify $B(aG)$ with a subalgebra of $\Mcb(A(G))$. 
\begin{proposition} \label{McbWAPprop2}
Let $G$ be a locally compact group. Then:
\begin{items}
\item $\Mcbo(A(G)) := \WAP_0(G) \cap \Mcb(A(G))$ is a closed ideal of $\Mcb(A(G))$;
\item $B(aG) = \AP(G) \cap \Mcb(A(G))$;
\item we have a direct sum decomposition
\begin{equation} \label{dsum}
  \Mcb(A(G)) = B(aG) \oplus \Mcbo(A(G)).
\end{equation}
\end{items}
\end{proposition}
\begin{proof}
(i) is obvious.
\par 
For (ii), the inclusion $B(aG) \subset \AP(G) \cap \Mcb(A(G))$ is clear. 
\par 
For the proof of the reverse inclusion, we make use of charaterizations of $\Mcb(A(G))$ from \cite{Spr}. Consider $\WAP(G) \tensor^{eh} \WAP(G)$ where $\tensor^{eh}$ stands for the extended Haagerup tensor product (see \cite{ER2} for the definition). We can view $\WAP(G) \tensor^{eh} \WAP(G)$ as a space of continuous functions on $G \times G$, and define 
\[
  \mathcal{VW}_\mathrm{inv}(G) := \{ u \in \WAP(G) \tensor^{eh} \WAP(G) : \text{$u(x,zy) = u(x,yz^{-1})$ for $x,y,z \in G$} \}.
\]
By \cite[Corollary 5.6]{Spr}, we have a completely isometric isomorphism between $\Mcb(A(G))$ and $\mathcal{VW}_\mathrm{inv}(G)$ via $\Mcb(A(G)) \ni \phi \to u_\phi$ with 
\[
  u_\phi(x,y) = \phi(xy^{-1}) \qquad (\phi \in \Mcb(A(G)), \, x,y \in G).
\]
Through the Gelfand transform of $\WAP(G)$ each $u \in \WAP(G) \tensor^{eh} \WAP(G)$ extends to hat $\hat{u} \in \mathcal{C}(wG) \tensor^{eh} \mathcal{C}(wG)$. It follows that 
\[
  \mathcal{VW}_\mathrm{inv}(G) = \{ u \in \WAP(G) \tensor^{eh} \WAP(G) : \text{$\hat{u}(sr,t) = \hat{u}(s,tr^\ast)$ for $s,t,r \in wG$} \}.
\]
\par 
Let $\phi \in \AP(G) \cap \Mcb(A(G))$, so that $\hat{u}_\phi \in \mathcal{VW}_\mathrm{inv}(G)$ with
\[
  \hat{\phi}(st^\ast) = \hat{u}_\phi(s,t) \qquad (s,t \in wG). 
\]
As $\phi \in \AP(G)$, we have 
\[
  \hat{\phi}(e_{\ker wG} \, s) = \hat{\phi}(s) = \hat{\phi}(s \, e_{\ker wG}) \qquad (s \in wG)
\]
(this follows from \cite[Theorem 2.20 and 2.22]{Bur} and the centrality of $e_{\ker wG}$). We conclude that 
\begin{equation} \label{APeq}
  \hat{u}_\phi(e_{\ker wG} \, s, e_{\ker wG} \, t) = \hat{\phi}(e_{\ker wG} \, s t^\ast \, e_{\ker wG}) = \hat{\phi}(st^\ast) = \hat{u}_\phi(s,t) \qquad (s,t \in wG).
\end{equation}
From the definition of $\tensor^{eh}$, we obtain $( \xi_i )_{i \in \mathbb{I}}$ and $( \eta_i )_{i \in \mathbb{I}}$ in $\WAP(G)$ for some index set $\mathbb{I}$ such that
$\hat{u}_\phi = \sum_{i \in \mathbb{I}} \hat{\xi}_i \tensor \hat{\eta}_i \in \WAP(G) \tensor^{eh} \WAP(G)$. From (\ref{APeq}) we conclude that 
$\hat{u}_\phi = \sum_{i \in \mathbb{I}} \hat{\xi}_i(e_{\ker wG} \, \cdot) \tensor \hat{\eta}_i(e_{\ker wG} \, \cdot ) \in \WAP(G) \tensor^{eh} \WAP(G)$. As $\hat{\xi}_i(e_{\ker wG} \, \cdot)$ and $\hat{\eta}_i(e_{\ker wG} \, \cdot )$ are almost periodic for $i \in \mathbb{I}$, we conclude that $\hat{u}_\phi \in \AP(G) \tensor^{eh} \AP(G) \cong \mathcal{C}(aG) \tensor^{eh} \mathcal{C}(aG)$
satisfying 
\[
  \hat{u}_\phi(xz,y) = \hat{u}_\phi(x,yz^{-1}) \qquad (x,y,z \in aG).
\]
It follows that $\phi \in \Mcb(A(aG)) = B(aG)$.
\par 
For (iii), let $\phi \in \Mcb(A(G))$. Then an inspection of the proof of (ii) yields that $\hat{\phi}(e_{\ker wG} \, \cdot ) \in B(aG)$, and obviously $\phi - \hat{\phi}(e_{\ker wG} \, \cdot ) \in \Mcbo(A(G))$.
\end{proof}
\begin{remarks}
\item The direct sum decomposition Proposition (\ref{dsum}) exists also in the category of operator spaces, i.e., the projections onto the summands are completely bounded.
\item Both $B(aG)$ and $\Mcbo(A(G))$ are dual Banach spaces. However, (\ref{dsum}) does not correspond to a decomposition of the predual of $\Mcb(A(G))$. Even though $\Mcbo(A(G))$ is weak$^\ast$ closed in $\Mcb(AG))$, $B(aG)$ rarely is. In fact, for many groups $B(aG)$ is weak$^\ast$ dense in $\Mcb(A(G))$: this fact is used in \cite{FRS} to establish the operator amenability of $\AMcb(A(\free_2))$.
\item In particular, the projection from $\Mcb(A(G))$ onto $B(aG)$ arising from (\ref{dsum}) is \emph{not} the restriction map $\Mcb(A(G)) \ni \phi \mapsto \hat{\phi} |_{\ker wG}$.
\end{remarks}
\par 
For the sake of notational simplicity, we shall from now on simply write $aG$ for $\ker wG$.
\section{Amenability of $\AMcb(G)$}
Let $G$ be a locally compact group such that $\AMcb(G)$ is amenable. By Proposition \ref{BrianVolker}, this means that $\chi_{\nabla(G)} \in \Mcb(A(G_d \times G_d))$. The direct sum decomposition (\ref{dsum}) yields a function in $B(a(G_d \times G_d) ) = B(a(G_d) \times a(G_d))$. Could this function be $\chi_{\nabla(a(G_d))}$ ? The answer is ``no''---except in trivial cases: if $G$ is a locally compact group such that $\chi_{\nabla(G)} \in \Cb(G \times G)$, then $G$ has to be discrete, which is very easy to see. More work is therefore required.
\par 
Given a locally compact group $G$, we denote its full group $\cstar$-algebra by $\cstar(G)$, and by $W^\ast(G)$ the enveloping von Neumann algebra, i.e., the second dual, of $\cstar(G)$. Let $\omega_G \!: G \to W^\ast(G)$ be the universal unitary representation of $G$. By the universal property of the full group $\cstar$-algebra, it induces a $^\ast$-homomorphism $\tilde{\omega}_{G_d} \!: \cstar^\ast(G_d) \to W^\ast(G)$. Consequently, $\tilde{\omega}_{G_d}^\ast$ maps $W^\ast(G)^\ast$ contractively into $\cstar(G_d)^\ast = B(G_d)$.
\par 
Dualizing (\ref{dsum}), we obtain a direct sum decomposition 
\[
  \Mcb(A(G))^\ast = W^\ast(aG) \oplus \Mcbo(A(G))^\ast.
\]
So, if $\phi \in M_{cb}(A(G))$, restricting it to $W^\ast(aG)$ yields a bounded linear functional on $W^\ast(aG)$, which is mapped by $\tilde{\omega}_{G_d}^\ast$ into $B((aG)_d)$. We define $\rho_{(aG)_d} \!: \Mcb(A(G)) \to B((aG)_d)$ by letting
\[
  \rho_{(aG)_d}(\phi)(x) := \langle \omega_{aG}(x), \phi \rangle \qquad (\phi \in \Mcb(A(G)), \, x \in G).
\]
\par 
We need three more technical lemmas:
\begin{lemma} \label{lem1}
Let $G$ be a locally compact group, and let $S \subset G$. Then 
\[
  \varcl{\iota_w(S)}^{wG} \cap aG = e_{aG} \, \varcl{\iota_w(S)}^{wG} = \varcl{\iota_a(S)}^{aG}
\]
holds.
\end{lemma}
\begin{proof}
As $\iota_a = p_{aG} \circ \iota_w$, it is clear that $e_{aG} \varcl{\iota_w(S)}^{wG} \subset \varcl{\iota_a(S)}^{aG}$. On the other hand, $e_{aG} \varcl{\iota_w(S)}^{wG}$ is a compact subset of $aG$ containing $\iota_a(S)$, so that $\varcl{\iota_a(S)}^{aG} \subset e_{aG} \varcl{\iota_w(S)}^{wG}$ as well. This proves the second equality.
\par
It is obvious that $\varcl{\iota_w(S)}^{wG} \cap aG \subset e_{aG} \varcl{\iota_w(S)}^{wG}$. For the reverse inclusion, let $s \in e_{aG} \, \varcl{\iota_w(S)}^{wG}$. Let $\mathcal{V}$ denote the set of neighborhoods of $s$ order by reversed set inclusion, i.e.,
\[
  V_1 \prec V_2 \defiff V_2 \subset V_1 \qquad (V_1, V_2 \in \mathcal{V}).
\]
Let $V \in \mathcal{V}$. As $s \in e_{aG} \, \varcl{\iota_w(S)}^{wG}$, there is $y_V \in G$ such that $e_{aG} \iota_w(y_v) \in V$. As multiplication in $wG$ is separately continuous, there is $x_V \in G$ such that $\iota_w(x_v y_V) = \iota_w(x_v) \iota_w(y_V) \in V$. It follows that $s = \lim_{V \in \mathcal{V}} \iota_w(x_v y_V)$, so that $s \in \varcl{\iota_w(S)}^{wG} \cap aG$. This proves the first equality.
\end{proof}
\begin{lemma} \label{lem2}
Let $G$ be a locally compact group, and let $S \subset G$ be such that $\chi_{S} \in \Mcb(A(G))$. Then $\rho_{(aG)_d}(\chi_S) = \chi_{\varcl{\iota_a(S)}^{aG}}$ holds, so that, in particular, $\chi_{\varcl{\iota_a(S)}^{aG}} \in B((aG)_d)$.
\end{lemma}
\begin{proof}
Let $x \in \chi_{\varcl{\iota_a(S)}^{aG}} = \varcl{\iota_w(S)}^{wG} \cap aG$. Then there is a net $( x_\alpha )_{\alpha \in \mathbb{A}}$ in $S$ such that $\omega_{aG}(x) = \lim_\alpha \iota_w(x_\alpha)$. It follows that 
\[
  \rho_{(aG)_d}(\chi_S)(x) = \lim_\alpha \chi_S(x_\alpha) = 1.
\]
\par
Conversely, let $x \in aG \setminus \varcl{\iota_w(S)}^{wG}$. Let $( x_\alpha )_{\alpha \in \mathbb{A}}$ be a net in $G$ such that $\omega_{aG}(x) = \lim_\alpha \iota_w(x_\alpha)$. Let $V$ be a neighborhood of $x$ such that $V \cap \varcl{\iota_w(S)}^{wG}= \void$. It follows that there is $\alpha_V \in \mathbb{A}$ such that $\iota_w(x_\alpha) \in V \subset wG \setminus \iota(S)$. It follows that 
\[
  \rho_{(aG)_d}(\chi_S)(x) = \lim_\alpha \chi_S(x_\alpha) = 1.
\]
In view of Lemma \ref{lem1}, this proves the claim.
\end{proof}
\begin{lemma} \label{lem3}
Let $G$ be a discrete group. Then $\varcl{\iota_a(\nabla(G)}^{aG \times aG} = \nabla(aG)$ holds.
\end{lemma}
\begin{proof}
It is straightforward to see that  $\varcl{\iota_a(\nabla(G)}^{aG \times aG} \subset \nabla(aG)$.
\par 
For the reverse inclusion, let $(x,x^{-1}) \in \nabla(aG)$. Let $(x_\alpha)_{\alpha \in \mathbb{A}}$ be such that $x = \lim_\alpha \iota_w(x_\alpha)$. It follows that 
\[
  x = \lim_\alpha \iota_w e_{aG} x_\alpha) \qquad\text{and}\qquad x = \lim_\alpha \iota_w(x_\alpha) e_{aG}.
\]
Continuity of inversion applied to the second limit yields $x^{-1} = \lim_\alpha e_{aG} \iota_w(x_\alpha^{-1})$, so that 
\[
  (x,e_{aG})  = \lim_\alpha e_{aG} (\iota_w(x_\alpha),e_{aG}) \qquad\text{and}\qquad (e_{aG}, x^{-1}) = \lim_\alpha (e_{aG}, e_{aG} \iota_w(x_\alpha^{-1})).
\]
As multiplication in $aG$ is jointly continuous, we obtain 
\[
  (x,x^{-1}) = \lim_\alpha \underbrace{(e_{aG},e_{aG})}_{=e_{aG \times aG}} (\iota_w(x_\alpha), \iota_w(x^{-1}_\alpha)).
\]
In view of Lemma \ref{lem1}, this means that $(x,x^{-1}) \in \varcl{\iota_a(\nabla(G)}^{aG \times aG}$.
\end{proof}
\par 
We are finally ready to prove the main result of this paper:
\begin{theorem} \label{mainthm}
Let $G$ be a locally compact group such that $\AMcb(G)$ is amenable. Then $a(G_d)$ is almost abelian. 
\end{theorem}
\begin{proof}
By Proposition \ref{BrianVolker}, we have $\chi_{\nabla(G)} \in \Mcb(A(G_d \times G_d))$, and from Lemmas \ref{lem2} and \ref{lem3}, it follows that $\chi_{\nabla(a(G_d)} \in B((a(G_d \times G_d))_d) = B((a(G_d))_d \times (a(G_d))_d)$.
\par 
The remainder of the argument is the same as in \cite{FR1}, \cite{RunPAMS}, or \cite{RU}; we sketch it here for the reader's convenience.
\par 
For $f \in A(a(G_d))$, define $\check{f} \in A(a(G_d))$ by letting 
\[
  \check{f}(x) := f(x^{-1}) \qquad (x \in G).
\]
The fact that $\chi_{\nabla(a(G_d)} \in B((a(G_d))_d \times (a(G_d))_d)$ implies that $A(a(G_d)) \ni f \mapsto \check{f}$ is completely bounded, which is possible only if $a(G_d)$ is almost abelian.
\end{proof}
\par 
The next two corollaries are immediate:
\begin{corollary} \label{cor1}
Let $G$ be a locally compact group such that $\AMcb(G)$ is amenable and that $G_d$ is maximally almost periodic. Then $G$ is almost abelian.
\end{corollary}
\begin{corollary} \label{cor2}
Let $G$ be a discrete, maximally almost periodic groups such that $\AMcb(G)$ is amenable. Then $G$ is almost abelian.
\end{corollary}
\begin{example}
The free group in two generators---$\free_2$---is maximally almost periodic, but not almost abelian. Therefore, $\AMcb(\free_2)$ is not amenable.
\end{example}
\par 
The example of $\free_2$ implies that \emph{if} there is a locally compact group $G$ with amenable $\AMcb(G)$, it will probably be very difficult to find:
\begin{corollary} \label{cor3}
Let $G$ be a locally compact group, such that $\AMcb(G)$ is amenable. Then one of the following holds:
\begin{alphitems}
\item $G$ is almost abelian;
\item $G$ is not amenable and does not have a closed subgroup isomorphic to $\free_2$.
\end{alphitems}
\end{corollary}
\begin{proof}
Suppose that $G$ is not almost abelian.
\par 
If $G$ is amenable, then $A(G) = \AMcb(G)$, so that $G$ is almost abelian by \cite{FR1}, which is a contradiction.
\par 
If $G$ contains $\free_2$ as a closed subgroup, then \cite[Theorem 1]{For} and basic hereditary properties of amenability (see \cite[Section 2.3]{RunAM}) imply that $\AMcb(\free_2)$ is amenable, which is also a contradiction.
\end{proof}
\begin{remarks}
\item The question if groups as in Corollary \ref{cor3}(b) actually do exist was already raised by J.\ von Neumann. A first example of such a group was eventually constructed by A.\ Yu.\ Ol'shanski\u{\i} (\cite{Olsh}). For references to further examples, see the Notes and Comments sections of \cite[Chapter 1]{RunAM}.
\item Let $G$ be a localdly compact group. For every multiplier $\phi$ of $A(G)$, the operator $M_\phi \!: A(G) \to A(G)$ is bounded by the closed graph theorem. Defining $\| \phi \|_M := \| M_\phi \|$ turns the multipliers of $A(G)$ into a Banach algebra denoted by $M(A(G))$. The closure of $A(G)$ in $M(A(G))$ is denoted by $A_M(G)$. We do not know if $A_M(\free_2)$ is amenable.
\end{remarks}
\renewcommand{\baselinestretch}{1.0}
\renewcommand{\baselinestretch}{1.2}
\dated
\vfill
\begin{tabbing}
\textit{Author's address}: \= Department of Mathematical and Statistical Sciences \\
                           \> University of Alberta \\
                           \> Edmonton, Alberta \\
                           \> Canada T6G 2G1 \\[\medskipamount]
\textit{E-mail}:           \> \texttt{vrunde@ualberta.ca}\\[\medskipamount]
\textit{URL}:              \> \texttt{http://sites.ualberta.ca/$^\sim$runde/}
\end{tabbing}
\end{document}